\documentclass{amsart}

\makeatletter
\let\@wraptoccontribs\wraptoccontribs
\makeatother

\usepackage{amsmath}
\usepackage{amsxtra}
\usepackage{amscd}
\usepackage{amsthm}
\usepackage{amsfonts}
\usepackage{amssymb}
\usepackage{mathtools}
\usepackage{eucal}
\usepackage[all]{xy}
\usepackage{graphicx}
\usepackage{tikz-cd}
\usepackage{mathrsfs}
\usepackage{euler}
\usepackage{color}
\usepackage{longtable}
\usepackage{caption}
\usepackage{enumerate}
\usepackage{faktor}
\usepackage{stmaryrd}
\usepackage{rotating}
\usepackage[T1]{fontenc} % improved font encoding
\usepackage[utf8]{inputenc} % for better handling of non-ASCII characters
\usepackage{tikz}
\usepackage{tikz-cd}
\usetikzlibrary{matrix,arrows,decorations.pathmorphing}
\usepackage{hyperref}
\usepackage{float}
\usepackage{cleveref}

\theoremstyle{plain}
\newtheorem{theorem}{Theorem}
\newtheorem{lemma}[theorem]{Lemma}
\newtheorem{proposition}[theorem]{Proposition}
\newtheorem{corollary}[theorem]{Corollary}
\theoremstyle{definition}

\theoremstyle{remark}

\newtheorem{ex}[theorem]{Example}
\numberwithin{theorem}{section}
\numberwithin{equation}{section}

\DeclareMathAlphabet{\mathcal}{OMS}{cmsy}{m}{n}
\DeclareFontShape{OT1}{cmr}{b}{up}{<-> ssub * cmr/b/n}{}
\DeclareFontShape{OT1}{cmr}{m}{up}{<-> ssub * cmr/m/n}{}

\newcommand{\C}{\mathbb{C}}

\newcommand{\F}{\mathbb{F}}

\newcommand{\Q}{\mathbb{Q}}

\newcommand{\Z}{\mathbb{Z}}

\newcommand{\calh}{\mathcal{H}}

\newcommand{\Ind}{\mathrm{Ind}}
\newcommand{\Rad}{\mathrm{Rad}}
\newcommand{\ov}[1]{{\overline{#1}}}
\newcommand{\onto}{\twoheadrightarrow}
\newcommand{\into}{\hookrightarrow}
\newcommand{\cO}{\mathcal{O}}
\DeclareMathOperator{\GL}{GL}

\newcommand{\Hom}{\mathrm{Hom}}

\newcommand{\Res}{\mathrm{Res}}
\newcommand{\End}{\mathrm{End}}

\newcommand{\soc}{\mathrm{soc}}
\newcommand{\cosoc}{\mathrm{cosoc}}

\newcommand{\triv}{\mathrm{triv}}

\newcommand{\Gal}{\mathrm{Gal}}

\newcommand{\nm}{\mathrm{Nm}}
\newcommand{\p}{\mathfrak{p}}

\newcommand{\paren}[1]{\mathopen{}\left(#1\right)\mathclose{}}
\newcommand{\set}[1]{\mathopen{}\left\{#1\right\}\mathclose{}}

\newcommand{\abrac}[1]{\mathopen{}\left\langle#1\right\rangle\mathclose{}}

\newcommand\restr[2]{{% we make the whole thing an ordinary symbol
  \left.\kern-\nulldelimiterspace % automatically resize the bar with \right
  #1 % the function
  \right|_{#2} % this is the delimiter
  }}
\newcommand{\Mid}{\,\middle|\,}

\title[Corrigendum]
{Corrigendum to ``Modular Gelfand pairs and multiplicity-free representations''}

\author{Patrick B. Allen}
\address[Patrick Allen]{Department of Mathematics and Statistics, McGill University}
\email{patrick.allen@mcgill.ca}

\author{Preston Wake}
\address[Preston Wake]{Department of Mathematics, Michigan State University}
\email{wakepres@msu.edu}

\author{Robin Zhang}
\address[Robin Zhang]{Department of Mathematics, Massachusetts Institute of Technology}
\email{robinz@mit.edu, robinzmit@gmail.com}

\date{May 6, 2026}

\begin{document}

\begin{abstract}
Some of the multiplicity-freeness results in ``Modular Gelfand pairs and multiplicity-free representations'' are stated in overly broad generality. We provide counterexamples and partial corrections.
\end{abstract}

% \keywords{Gelfand pairs, multiplicity-free triples, multiplicity one, modular representations, Hecke algebras, self-projective modules, self-injective modules}
% \subjclass{Primary: 20C20; Secondary: 11F70, 16D40, 16D50, 20C08, 20G05, 20G40.}

\maketitle

\setcounter{tocdepth}{1}

\tableofcontents

%%%%%%%%%%%%%%%%%%%%%%%%%%%% Introduction %%%%%%%%%%%%%%%%%%%%%%%%%%%%

\section{Introduction}

The purpose of this article is to correct some errors in \cite{zhang2024}. One of the main results of that paper is \cite[Corollary 4.8]{zhang2024}, which claims that for a finite group $G$ with a subgroup $H$ and a complex representation $\rho$ of $H$ such that the mod-$p$ reduction $\overline{\rho}$ is irreducible, if $\Ind_H^G(\rho)$ is multiplicity-free, then the socle and cosocle of $\Ind_H^G(\overline{\rho})$ are multiplicity-free. This is false in general, as the following example (which is explained fully in Section \ref{subsec:examples} below) illustrates.
\begin{ex}
\label{ex:dihedral intro}
    Let $G=D_{30}$ be the dihedral group of order $30$, $H$ be the Sylow-$5$ subgroup of $H$, $\rho:H \to \C^\times$ be a non-trivial character, and let $p=5$. Then $\Ind_H^G(\rho)$ is multiplicity-free, but both the socle and cosocle of $\Ind_H^G (\overline{\rho})$ are not.
\end{ex}

The strategy of \cite{zhang2024} is to extend the classical observation that a semisimple module $M$ is multiplicity-free if and only if $\End(M)$ is commutative. Neither implication is true without the semisimplicity assumption (as \cite[Non-examples 3.4, 3.5]{zhang2024} illustrate). However, if the natural maps $\End(M) \to \End(\soc(M))$ and $\End(M) \to \End(\cosoc(M))$ are surjective, then this does show that the commutativity of $\End(M)$ implies that the semisimple modules $\soc(M)$ and $\cosoc(M)$ are multiplicity-free. These two respective maps are surjective if $M$ is injective or projective, respectively, or even if $M$ satisfies the weaker condition of being \emph{self-injective} or \emph{self-projective}, respectively. The proof of \cite[Corollary 4.8]{zhang2024} proceeds by showing that, if $\Ind_H^G(\rho)$ has commutative endomorphism ring, then so does $\Ind_H^G(\overline{\rho})$ (in the proof of \cite[Corollary 4.8]{zhang2024}), and, if $\overline{\rho}$ is irreducible, then $\Ind_H^G(\overline{\rho})$ is both self-injective and self-projective (\cite[Lemma 4.3(iii-iv)]{zhang2024}). However, the proofs of \cite[Corollary 4.8 and Lemma 4.3(iii-iv)]{zhang2024} contain errors:
\begin{itemize}
    \item In the proof of [Zha24, Corollary 4.8], the claim that $\End_G(\Ind_H^G(\overline{\rho}))$ is commutative if $\End_G(\Ind_H^G(\rho))$ has commutative endomorphism ring is deduced from an erroneous assertion that there is always a surjection of Hecke algebras $\calh(G, H, \rho, \C) \twoheadrightarrow \calh(G, H, \overline{\rho}, \overline{\F_p})$.
    \item In the proof of [Zha24, Lemma 4.3(iii-iv)], the claim that $\Ind_H^G(\overline{\rho})$ is self-projective (resp. self-injective) if $\overline{\rho}$ is self-projective (resp. self-injective) is deduced from an erroneous application of Mackey's restriction formula.
\end{itemize}
Our corrected versions can be described as follows:
\begin{itemize}
    \item If $\Ind_H^G(\rho)$ is multiplicity free, then $\End_G(\Ind_H^G(\overline{\rho}))$ need not be commutative in general, but we give a necessary and sufficient criterion for this to hold (Proposition \ref{prop:reduction of end ind}). We give examples, including Example \ref{ex:dihedral intro}, where this criterion is not satisfied.
    \item If $\overline{\rho}$ is irreducible, then $\Ind_H^G(\overline{\rho})$ need not be self-projective or self-injective. We prove a necessary criterion (Lemma \ref{lem:necessary for self-proj}) and use it exhibit counterexmples to \cite[Lemma 4.3(iii-iv)]{zhang2024}. We give simple sufficient criteria (Lemma \ref{lem:selfproj criteria}) that include:
    \begin{itemize}
        \item $p \nmid \# H$, or
        \item $H$ is normal in $G$,
    \end{itemize}
    and also review a criterion due to Vign\'{e}ras. Our necessary criterion is very far from these sufficient criteria and the problem of precisely determining when $\Ind_H^G(\overline{\rho})$ is self-projective or self-injective seems difficult.
\end{itemize}
When the sufficient criteria for these two points are both met, then one can conclude that the socle and cosocle of $\Ind_H^G(\overline{\rho})$ are multiplicity free. One simple case where the criteria are both met is when $H$ has prime-to-$p$ order and $\rho$ is a character. This case recovers \cite[Corollaries 5.2 and 5.4]{zhang2024} about multiplicity-freeness of modular Gelfand--Graev modules (indeed, as noted in \cite[Remark 5.3]{zhang2024}, these statements are vacuous in the equal-characteristic case $\ell \mid q$, while the relevant group $H$ has prime-to-$p$ order in the unequal characteristic case $\ell \nmid q$). This also recovers the applications of \cite{zhang2024} used in \cite[Section 4.1]{gamma}.

% \subsection{Further comments} 
To avoid possible confusion about how errors in \cite{zhang2024} affect the results here, we have made this article completely independent and do not assume that the reader is familiar with \cite{zhang2024}. We have attempted to reproduce corrected versions of the most important results of \cite{zhang2024} here, so that the reader interested in those results need only look at this paper (although we note that \cite[Sections 2, 3, and 6]{zhang2024} are not affected by the errors that we identify).
For the most part, we only deal with finite groups and
have not attempted to extend the results to infinite compact groups,
which are the subject of \cite[Corollary 1.10 and Sections 4.3, 5.2]{zhang2024}.
However, in Section \ref{sec:quaternion},
we reproduce a proof for the $\ell=p$ case of \cite[Corollary 1.10]{zhang2024}
that we learned from Dipendra Prasad.
We do not know whether the $\ell \neq p$ case of \cite[Corollary 1.10]{zhang2024} is true
because the methods we develop here are not sufficient
to determine whether self-injectivity holds.

\subsection*{Acknowledgements}
R.~Z.~is very grateful to his coauthors for identifying and communicating the errors
and for agreeing to work together on this note.
He thanks Karol Koziol for identifying some of the errors and pointing out Example \ref{ex:gl2-unipotent}; Ofer Gabber for explaining the equivalence of self-projectivity and almost-projectivity when there are enough projectives;
and Dipendra Prasad for pointing out how to recover Corollary \ref{cor:qda}. We thank the anonymous referee for their careful reading and comments.

P. A. acknowledges the support of the Natural Sciences and Engineering Research Council of Canada (NSERC), [funding reference number RGPIN-2020-05915].
P. A. a été financée par le Conseil de recherches en sciences naturelles et en génie du Canada (CRSNG), [numéro de référence RGPIN-2020-05915].
P.~W.~is supported by National Science Foundation CAREER Grant No.~DMS-2337830.
R.~Z.~is supported by National Science Foundation Grant No.~DMS-2303280.

\section{Endomorphisms of induced representations}
We use Mackey theory to work out when the reduction map on twisted Hecke algebras is surjective.

\subsection{Reduction modulo $p$ on endomorphisms of inductions} 
Let $p$ be a prime number, and let $K$ be a finite extension of $\Q_p$ with valuation ring $\cO$, residue field $k$, and uniformizer $\varpi$. Whenever we discuss $K$-valued representations of a finite group $G$, we assume that $K$ is large enough that all the irreducible representations of $G$ are defined over $K$.
For an $\cO$-module $M$, write $\ov{M}=M/\varpi M$ and $M_K=M\otimes_\cO K$.

\begin{lemma}
\label{lem:mod p end}
Let $R$ be an $\cO$-algebra that is finitely generated and free as an $\cO$-module, and let $M$ be a finitely generated $R$-module that is free as an $\cO$-module. Then the map
\[
\ov{\End_R(M)} \to \End_{\ov{R}}(\ov{M})
\]
induced by $\phi \mapsto \ov{\phi}$ is injective. It is surjective if and only if
\[
\dim_K(\End_{R_K}(M_K)) \ge \dim_k(\End_{\ov{R}}(\ov{M})).
\]
\end{lemma}
\begin{proof}
Applying the functor $\Hom_R(M,-)$ to the exact sequence
    \[
    0 \to M \xrightarrow{\varpi} M \to \ov{M} \to 0
    \]
yields an exact sequence
\[
0 \to \End_R(M) \xrightarrow{\varpi} \End_R(M) \to \Hom_R(M,\ov{M}).
\]
Since $\Hom_R(M,\ov{M}) = \End_{\ov{R}}(\ov{M})$, this proves the first statement, and also shows that $\End_R(M)$ is a free $\cO$-module. Since the map $\cO \to K$ is flat, the map
\[
\End_R(M)_K \to \End_{R_K}(M_K)
\]
is an isomorphism, so
$\End_R(M)$ is a free $\cO$-module of rank $\dim_K(\End_{R_K}(M_K))$. This implies that $\dim_k(\ov{\End_R(M)})=\dim_K(\End_{R_K}(M_K))$, which, together with the first statement, proves the second statement.
\end{proof}

Let $G$ be a finite group, $H$ be a subgroup of $G$, and $\rho: H \to \GL_n(L)$ be a finite-dimensional representation over a field $L$. For an element $g \in H$, let $H_g = H \cap g^{-1}Hg$, and let $\rho^g$ be the representation of $H_g$ given by $\rho^g(x)=\rho(gx g^{-1})$.  Let $[H \backslash G /H] \subset G$ be a set of coset representatives. Then Mackey's formula states that
\[
    \Res_H^G(\Ind_H^G(\rho)) \cong \bigoplus_{s \in [H \backslash G /H]} \Ind_{H_s}^H(\rho^s).
\]
Applying Frobenius reciprocity twice then implies that
\begin{equation}
\label{eq:mackey}
    \End_{L[G]}(\Ind_H^G(\rho)) = \bigoplus_{s \in [H \backslash G /H]} \Hom_{L[H_s]}(\rho,\rho^s)
\end{equation}

\begin{proposition}
\label{prop:reduction of end ind}
Let $\rho: H \to \GL_n(\cO)$ be an $\cO$-lattice in a representation $\rho_K$ over $K$. Then the map
\[
    \ov{\End_{\cO[G]}(\Ind_H^G(\rho))} \to \End_{k[G]}(\Ind_H^G(\ov{\rho}))
\]
is injective and is an isomorphism if and only if
\begin{equation}
\label{eq:end ind surj inequality}
    \sum_{s \in [H \backslash G /H]} \dim_K(\Hom_{K[H_s]}(\rho_K,\rho^s_K)) \ge \sum_{s \in [H \backslash G /H]} \dim_k(\Hom_{k[H_s]}(\bar{\rho},\bar{\rho}^s)).
\end{equation}
Moreover, if $n=1$, then \eqref{eq:end ind surj inequality} holds if and only if $\Res_{H_s}^H(\rho_K)^{-1} \, \rho_K^s$ is either trivial or does not have $p$-power order for all $s \in [H \backslash G /H]$. In particular, if $n=1$ and $\rho$ has prime-to-$p$ order, then \eqref{eq:end ind surj inequality} holds.
\end{proposition}
\begin{proof}
    The first statement is simply a combination of Lemma \ref{lem:mod p end} and Mackey's formula \eqref{eq:mackey}. For $n=1$, \eqref{eq:end ind surj inequality} simplifies to
    \[
    \#\set{s \in [H \backslash G /H] \ \Mid \ \Res_{H_s}^H(\rho_K) \cong \rho_K^s} \geq \#\set{s \in [H \backslash G /H] \ \Mid \ \Res_{H_s}^H(\ov{\rho}) \cong \ov{\rho}^s}.
    \]
This can fail only if there is a $s \in [H \backslash G /H]$ such that $\Res_{H_s}^H(\rho_K)^{-1} \, \rho_K^s$ is non-trivial and its reduction $\Res_{H_s}^H(\ov{\rho})^{-1} \,\ov{\rho}^s$ is trivial. This only occurs if and only if $\Res_{H_s}^H(\rho_K)^{-1} \, \rho_K^s$ is non-trivial and has $p$-power-order. Finally, since the order of $\Res_{H_s}^H(\rho_K)^{-1} \, \rho_K^s$ divides the order of $\rho$, the final statement follows.
\end{proof}

\begin{corollary}
\label{cor:commutativity of end ind}
Let $\rho: H \to \GL_n(\cO)$ be an $\cO$-lattice in a representation $\rho_K$ over $K$, and suppose that $\Ind_H^G(\rho_K)$ is multiplicity-free and that \eqref{eq:end ind surj inequality} holds. Then $\End_{k[G]}(\Ind_H^G(\overline{\rho}))$ is commutative.
\end{corollary}
\begin{proof}
Since $\Ind_H^G(\rho_K)$ is multiplicity-free, $\End_{K[G]}(\Ind_H^G(\rho_K))$ is commutative. Hence the subquotient ring $\ov{\End_{\cO[G]}(\Ind_H^G(\rho))}$ is a also commutative. By Proposition \ref{prop:reduction of end ind},  $\End_{k[G]}(\Ind_H^G(\overline{\rho}))$ is isomorphic to $\ov{\End_{\cO[G]}(\Ind_H^G(\rho))}$, so it is also commutative.
\end{proof}

\subsection{Examples}
\label{subsec:examples}
We give some examples where \eqref{eq:end ind surj inequality} fails and where the mod-$p$ twisted Hecke algebra is noncommutative even where the characteristic zero version is commutative.

\begin{ex}
    \label{ex:gl2-unipotent}
    Let $p$ be an odd prime, $\cO=\Z_p[\zeta_p]$, $G=\GL_2(\F_p)$, $H$ be the subgroup of upper-triangular unipotent matrices, and $\rho:H \to \cO^\times$ be a nontrivial character. By \eqref{eq:mackey}, the dimension of $\End_{K[G]}(\Ind_H^G(\rho))$ is equal to 
    \[
    \#\{s \in [H \backslash G/H] \ \mid \ s \not\in N_G(H) \text{ or } s\in C_G(H)\}. 
    \]
    Using the union of the diagonal matrices and the anti-diagonal matrices as the choice of $[H \backslash G/H]$, this dimension is easily seen to be $(p-1)^2+(p-1)$.
    But, since $\rho$ has order $p$, the reduction $\overline{\rho}$ is trivial, so  $\End_{k[G]}(\Ind_H^G(\triv_H)) = k[H \backslash G/H]$ has dimension $2(p-1)^2$. Hence this is an example where \eqref{eq:end ind surj inequality} fails.
\end{ex}

We now come to Example \ref{ex:dihedral intro} from the introduction. We first give a simple criterion for $\End_{K[G]}(\Ind_H^G(\rho))$ to be commutative.

\begin{lemma}
    Let $N$ be a normal subgroup of $G$ and let $\rho: N \to K^\times$ be a character. Since $N$ is normal, $\rho^s$ is also a character of $N$ for all $s \in G$. Let $H=\{s \in G \ | \ \rho^s=\rho\}$ be the stabilizer of $\rho$. If $[H:N] <4$, then $\End_{K[G]}\Ind_N^G(\rho)$ is commutative (and $\Ind_N^G(\rho)$ is multiplicity-free).
\end{lemma}
\begin{proof}
    Since $N$ is normal, $N \backslash G/N=G/N$. By \eqref{eq:mackey}, 
    \[
    \dim_K\End_{K[G]}\Ind_N^G(\rho) =\# \{s \in G/N \ | \ \rho^s=\rho\} = [H:N].
    \]
    Since $\End_{K[G]}\Ind_N^G(\rho)$ is a semisimple $K$-algebra, the fact that this dimension is less than $4$ implies that it is commutative.
\end{proof}

We also have a simple criterion for the socle and cosocle of $\Ind_H^G(\overline{\rho})$ to have multiplicity.

\begin{lemma}
    Let $G$ be a group with a normal $p$-Sylow subgroup $N$ such that $G/N$ is non-abelian. Then for every character $\rho:N \to K^\times$, the induction $\Ind_N^G(\overline{\rho})$ is semisimple and is not multiplicity free.
\end{lemma}
\begin{proof}
    Since $N$ is a $p$-group, $\overline{\rho}=\triv_N$ is the trivial character.
    Since $N$ is normal, $\Ind_N^G(\triv_N)$ is the inflation of the regular representation $k[G/N]$. Because $p \nmid \#(G/N)$, the usual formula for the regular representation in characteristic zero,
    \[
    k[G/N] \cong \bigoplus_{\pi \in \mathrm{Irr}(G/N)} \pi^{\oplus \dim_k(\pi)},
    \]
    holds over $k$. Since $G/N$ is non-abelian, it has at least one irreducible representation that is not one-dimensional, so $k[G/N]$ is not multiplicity free.
\end{proof}

Combining these two lemmas gives the following.

\begin{corollary}
    \label{cor:sylow}
    Let $G$ be a group that has a normal $p$-Sylow subgroup $N$ such that $G/N$ is non-abelian. Let $\rho:N \to K^\times$ be a character and let $H$ be the stabilizer of $\rho$. If $[H:N]<4$, then $\Ind_N^G(\rho)$ is multiplicity-free but $\Ind_N^G(\overline{\rho})$ is not.
\end{corollary}

A minimal example of a group satisfying the hypothesis of Corollary \ref{cor:sylow} is $G=D_{30}$, the dihedral group of order $30$, with $p=5$ as in Example \ref{ex:dihedral intro}.
The $5$-Sylow subgroup $N$ is normal and for every nontrivial $\rho:N \to K^\times$, the stabilizer $H$ of $\rho$ is the centralizer $C_G(N)$ of $N$, which is cyclic of order $15$. Hence $[H:N]=3<4$, and Corollary \ref{cor:sylow} applies.

\section{Self-projectivity of induced representations}

In this section, we consider the question of when an induced representation is a self-projective (or self-injective) module. Self-projectivity is a weak notion of projectivity that is enough to ensure that the map
\[
\End(M) \to \Hom(M,\cosoc(M)) \cong \End(\cosoc(M)) 
\]
is surjective.

For this section, $k$ denotes a field and $R$ a finite-dimensional $k$-algebra (unital but not necessarily commutative). A typical example is $R=k[G]$ for a finite group $G$. All $R$-modules considered are assumed to be finitely-generated left $R$-modules.

\subsection{Projective covers and injective envelopes}
Recall that an $R$-module surjection $\pi: N \onto M$ is called \emph{essential} if there is no proper submodule $Q \subsetneq N$ such that $\pi(Q)=M$. Equivalently, $\pi$ is essential if, for every proper submodule $Q \subsetneq N$, we have $Q+\ker(\pi) \subsetneq N$. An essential surjection $\pi: P \to M$ where $P$ is a projective module is called a \emph{projective cover} (we also abuse notation slightly and refer to the projective module $P$ as a projective cover of $M$).

\begin{lemma}
\label{lem:proj cover no homs}
    Let $M$ be a simple $R$-module and $P$ be a projective cover of $M$. If $N$ is a finitely generated $R$-module such that no subquotient of $N$ is isomorphic to $M$, then $\Hom_R(P,N)=0$.
\end{lemma}
\begin{proof}
    Note that the cosocle $P/\Rad(P)$ of $P$ is isomorphic to $M$. First consider the case that $N$ is semisimple; then the assumption is that $\Hom_R(M,N)=0$. But any homomorphism $P \to N$ factors through the cosocle $M$ of $P$, so it must be zero. Now consider general $N$ and let $f: P \to N$ be a homomorphism. By the semisimple case, the composition $P \xrightarrow{f} N \to N/\Rad(R)N$ is zero. Hence the image of $f$ is contained in $\Rad(R)N$, and the result follows by induction on the length of the radical series of $N$.
\end{proof}

Dually, a submodule $N \subseteq M$ is called \emph{essential} if there is no non-zero submodule $Q \subseteq M$ such that $N \cap Q=0$. If $I$ is an injective $R$-module such that $M$ is an essential submodule of $I$, then $I$ is called an \emph{injective envelope} of $M$.

\subsection{Self-projectivity and self-injectivity}
Let $M$ and $N$ be $R$-modules. Recall that an $R$-module $M$ is called \emph{$N$-projective} (resp.~\emph{$N$-injective}) if the functor $\Hom_R(M,-)$ (resp.~$\Hom_R(-,M)$) sends every short exact sequence of the form
\[
0 \to N' \to N \to N'' \to 0
\]
to a short exact sequence. Of course, $M$ is projective if and only if $M$ is $N$-projective for all $N$, and similarly for injectivity. Similarly, $N$ is semisimple if and only if every module $M$ is $N$-projective if and only if every module $M$ is $N$-injective. A module $M$ is called \emph{self-projective} (resp.~\emph{self-injective}) if $M$ is $M$-projective (resp.~$M$-injective).

In the following lemma, the implication $(2) \Rightarrow(3)$ was proved by Arabia \cite[Appendix, Proposition 7]{arabia-vigneras} (where property (3) is termed ``quasi-projectif'') in an arbitrary abelian category. We learned from Ofer Gabber that the converse holds for modules that have a projective cover.

\begin{lemma}
\label{lem:selfproj = almost proj}
Let $M$ be a finitely generated $R$-module. The following are equivalent:
\begin{enumerate}
    \item There is a projective cover $P$ of $M$ such that 
    \[
    \dim_k\End_R(M) =\dim_k\Hom_R(P,M).
    \]
    \item There is a projective $R$-module $P$ and a surjective homomorphism $\pi: P \onto M$ such that the induced injective map
    \[
    \End_R(M)\to \Hom_R(P,M) 
    \]
    is an isomorphism.
    \item $M$ is self-projective.
\end{enumerate}
\end{lemma}
\begin{proof}
    Clearly (1) implies (2), and $(2)$ implies $(3)$ by \cite[Appendix, Proposition 7]{arabia-vigneras}. Assume (3) and let $\pi:P \onto M$ be a projective cover with $K=\ker(\pi)$.  We will show that every $R$-module homomorphism $f:P \to M$ factors through $\pi$, which implies that the injective map $\End_R(M)\to \Hom_R(P,M)$ induced by $\pi$ is surjective, proving (1). Let $f:P \to M$ be a homomorphism, let $j:M \onto M/f(K)$ be the quotient map, and note that the kernel of the composition $j \circ f$ 
    contains $K$, so it can be written as $j \circ f =\ov{f} \circ\pi$ for a homomorphism $\ov{f}: M \to M/f(K)$.
    By (3), there is a map $g : M \to M$ such that $\ov{f} = j \circ g$. Together with the previous equation, this implies $j \circ f = j \circ g \circ \pi$. Letting $\tilde{g}=g\circ \pi$, this implies that the image of $f-\tilde{g} :P \to M$ is contained in $f(K)$. But the map $f-\tilde g$ agrees with $f$ on $K$, so the image of $f-\tilde{g}$ equals $f(K)$. This implies $P=\ker(f-\tilde g)+K$, so, by the definition of essential, $P=\ker(f-\tilde g)$. Hence $f=\tilde{g}$ and $f$ factors through $\pi$. 
\end{proof}

There is a dual lemma about self-injectivity, whose proof is almost the exact dual of the proof of Lemma \ref{lem:selfproj = almost proj}.

\begin{lemma}
\label{lem:almost injective and self-injective}
    Let $M$ be a finitely generated $R$-module. Then $M$ is self-injective if and only if there is an injective envelope $I$ of $M$ such that $\dim_k\End_R(M) =\dim_k\Hom_R(M,I).$

\end{lemma}

\subsection{Failure of self-projectivity and self-injectivity of certain inductions}
In this section, let $p$ be a prime and let $k$ be a field of characteristic $p$.
\begin{lemma}
\label{lem:necessary for self-proj}
    Let $G$ be a finite group and let $H$ be a non-normal subgroup of $G$. Let $M= \Ind_H^G(\triv_H)$ be the induction of the trivial representation of $H$. Suppose that $k[G]$ is a projective cover (resp.~injective envelope) of $M$. Then $M$ is not self-projective (resp.~self-injective). 
\end{lemma}
\begin{proof}
    Suppose that $k[G]$ is a projective cover of $M$. By Lemma \ref{lem:selfproj = almost proj}, to show that $M$ is not self-projective, it is enough to show that 
    \[
    \dim_k\End_{k[G]}(M) < \dim_k \Hom_{k[G]}(k[G],M)
    \]
    But $\Hom_{k[G]}(k[G],M)=M$ has dimension $[G:H]$ and $\End_{k[G]}(M)=k[H \backslash G/H]$ has dimension $\#(H \backslash G/H)$, which is strictly smaller than $[G:H]$ by the non-normality assumption.

    Now suppose that $k[G]$ is an injective envelope of $M$. By Lemma \ref{lem:almost injective and self-injective}, to show that $M$ is not self-injective, it is enough to show that 
    \[
    \dim_k\End_{k[G]}(M) < \dim_k \Hom_{k[G]}(M,k[G]).
    \]
    But $k[G]$ is self-dual, so $\Hom_{k[G]}(M,k[G]) \cong \Hom_{k[G]}(k[G],M^*)=M^*$ and the same arguments as in the projective case give the desired inequality.
\end{proof}

\begin{proposition}
If $G$ is a $p$-group and $H$ is a non-normal subgroup of $G$, then $\Ind_H^G(\triv_H)$ is neither self-projective nor self-injective.
\end{proposition}
\begin{proof}
    By the previous lemma, it is enough to show that $k[G]$ is a projective cover and injective envelope of $\Ind_H^G(\triv_H)$. Inducing the surjection $k[H] \onto k$ shows that $k[G]$ surjects onto $\Ind_H^G(\triv_H)$, so the projective cover of $\Ind_H^G(\triv_H)$ is a direct summand of $k[G]$. Similarly, inducing the injective $k \into k[H]$ shows that $\Ind_H^G(\triv_H)$ can be embedded in $k[G]$, so the injective envelope of $\Ind_H^G(\triv_H)$ is a direct summand of $k[G]$. But, since $k[G]$ is local, it has no proper direct summands.
\end{proof}

\subsection{Criteria for self-projectivity of inductions}
\label{subsec:vigneras}
Although the last section shows that inductions of simple modules are not self-projective in general, there are some simple criteria for when they are. We continue to assume that $p$ is a prime and $k$ is a field of characteristic $p$.

\begin{lemma}
\label{lem:selfproj criteria}
    Let $G$ be a finite group with a subgroup $H$ and let $M$ be a $k[H]$-module. Then $\Ind_H^G(M)$ is self-projective if any one of the following holds:
    \begin{enumerate}
        \item $M$ is $\Res_H^G(\Ind_H^G(M))$-projective.
        \item $p \nmid \#H$.
        \item $M$ is semisimple and $H$ is normal in $G$.
    \end{enumerate}
\end{lemma}
\begin{proof}
If (1) holds, then $\Ind_H^G(M)$ is self-projective by a Frobenius-reciprocity argument (see \cite[Corollary 4.2(ii)]{zhang2024}). If (2) holds, then $k[H]$ is a semisimple ring, so $M$ is a projective $k[H]$-module and (1) holds. If (3) holds, then $\Res_H^G(\Ind_H^G(M))$ is a semisimple $k[H]$-module by Mackey's formula, so every $k[H]$-module $N$ is $\Res_H^G(\Ind_H^G(M))$-projective and (1) holds again.
\end{proof}

The following criterion is a slight generalization of \cite[Lemma 3.1]{vigneras-2001}, in that we show that the assumption there that $M_\rho$ is a direct factor of $\Res_H^G(M)$ is unnecessary.
\begin{lemma}
    Let $G$ be a finite group with a subgroup $H$. Let $\rho$ be a simple $k[H]$-module and let $M=\Ind_H^G(\rho)$. Let $M_\rho$ be the largest $k[H]$-submodule of $\Res_H^G(M)$ that is isomorphic to a direct sum of copies of $\rho$, and let $M^\rho=\Res_H^G(M)/M_\rho$. If no subquotient of $M^\rho$ is isomorphic to $\rho$, then $M$ is self-projective.
\end{lemma}
\begin{proof}
    Let $P_\rho$ be a projective cover of $\rho$, and let $P=\Ind_H^G(P_\rho)$. Then $P$ is a projective $k[G]$-module and there is a surjection $P \onto M$. By Frobenius reciprocity, 
    \[
    \End_{k[G]}(M)=\Hom_{k[H]}(\rho, \Res_H^G(M))= \Hom_{k[H]}(\rho,M_\rho)
    \]
    and 
    \[
    \Hom_{k[G]}(P,M) = \Hom_{k[H]}(P_\rho, \Res_H^G(M)).
    \]
    Since $\Hom_{k[H]}(P_\rho,M^\rho)=0$ by Lemma \ref{lem:proj cover no homs}, it follows that
    \[
    \Hom_{k[H]}(P_\rho, \Res_H^G(M))=\Hom_{k[H]}(P_\rho,M_\rho).
    \]
    But every homomorphism $P_\rho \to M_\rho$ factors through the cosocle $\rho$ of $P_\rho$, so $\Hom_{k[H]}(P_\rho,M_\rho) \cong \Hom_{k[H]}(\rho,M_\rho)$. Altogether, this shows that $\End_{k[G]}(M)$ and $\Hom_{k[G]}(P,M)$ are isomorphic, so $M$ is self-projective by Lemma \ref{lem:selfproj = almost proj}
\end{proof}

Vign\'{e}ras uses this criterion to show \cite[Corollary 5.2]{vigneras-2001} that $\Ind_P^G(\sigma)$ is self-projective if $G=\mathbf{G}(F)$ is the $F$-rational points of a reductive group for a finite field $F$ with $\mathrm{char}(F)\ne \mathrm{char}(k)$, $P=MU$ is a parabolic subgroup with Levi $M$ and unipotent radical $U$, and $\sigma$ is a cuspidal representation of $M$.

%%%%%%%%%%%%%%%%%%%%%%%%%%%%%%%%%%%%%%%%%%%%%%%%%%%%%%%%%%%%%%%%%%%

\section{Quaternion division algebras}
\label{sec:quaternion}
Let $\ell$ and $p$ be two primes.
Let $F$ be a $p$-adic local field
with residue field $\F_{p^f}$
and let $D$ be the quaternion division algebra over $F$.
The main application of \cite{zhang2024} for infinite groups
is \cite[Corollary 5.6]{zhang2024}, which asserts the uniqueness of trilinear forms on irreducible smooth representations of $D^\times$
with coefficients in $k = \overline{\F}_\ell$. As a result of the errors in \cite{zhang2024} identified here, we do not know whether or not \cite[Corollary 5.6]{zhang2024} is true for $\ell \ne p$, but we learned from Dipendra Prasad that the $\ell=p$ case
can be recovered using results of Vign\'{e}ras \cite{vigneras-1989-quaternion}, as we now explain. For the remainder of the section, we take $\ell=p$ to be an odd prime and $k=\ov\F_p$.

\begin{corollary}[{\cite[Corollary 5.6]{zhang2024}} for odd $\ell = p$]
  \label{cor:qda}
  Let $p$ be an odd prime, $k = \overline{\F}_p$,
  $F$ be a $p$-adic local field, and $D$ be the quaternion division
  algebra over $F$.
  If $\pi_1, \pi_2, \pi_3$ are irreducible smooth
  $k$-representations of $D^\times$, then there exists at most one
  (up to isomorphism) non-zero $D^\times$-invariant linear form
  on $\pi_1 \otimes \pi_2 \otimes \pi_3$ over~$k$.
\end{corollary}

The proof uses the classification of irreducible smooth $k$-representations~of $D^\times$ due to Vign\'eras \cite{vigneras-1989-quaternion}, which we now recall.
Let $\cO_F$ be the ring of integers of $F$ and $\varpi_F$ be a uniformizer. Let $\cO_D$ be the ring of integers in $D$, let $\p_D \subset\cO_D$ be the maximal ideal, and let $U_D^1 := 1+\p_D \subset \cO_D^\times$. Since $U_D^1$ is a normal pro-$p$ subgroup of $D^\times$,  every irreducible smooth $k$-representation of $D^\times$ is inflated from the group $D^\times/U_D^1$ (see \cite[Corollaire 8.1, pg.259]{vigneras-1989-quaternion}). This is not a finite group, but we can arrange for such a representation to factor through the finite quotient $D^\times/\varpi_F^\Z U_D^1$ by twisting. Indeed, for an irreducible smooth $k$-representation $\pi$ of $D^\times$, let $\omega_\pi:F^\times \to k^\times$ be its central character; if $\chi: F^\times \to k^\times$ is a character such that $\chi(\varpi_F)^2=\omega_\pi(\varpi_F)$ and $\chi_{D^\times}:D^\times \to k^\times$ is the composition of $\chi$ with the reduced norm, then the representation $\pi \otimes \chi_{D^\times}^{-1}$ is trivial on $\varpi_F$, hence it factors through the finite group $D^\times/\varpi_F^\Z U_D^1$. The valuation map defines an exact sequence
\[
    1 \to \cO_D^\times/U_D^1 \to D^\times/\varpi_F^\Z U_D^1 \xrightarrow{\mathrm{val}} \Z/2\Z \to 1
\]
that splits (indeed, there is a generator $\varpi_D$ of $\p_D$ such that $\varpi_D^2=\varpi_F$), and we note that $\cO_D^\times/U_D^1$ is isomorphic to $\F_{p^{2f}}^\times$. Let $G=D^\times/\varpi_F^\Z U_D^1$ be this dihedral group, and write it as $G\cong C \rtimes \abrac{\sigma}$, where $C=\F_{p^{2f}}^\times$ and $\sigma$ acts by the Frobenius over $\F_{p^f}$. Vign\'{e}ras classifies irreducible smooth $k$-representations of $D^\times$ using the fact that (up to twist) every such representation is inflated from $G$.

To state the classification precisely, we need more notation. Let $E/F$ be the quadratic unramified extension and fix an $F$-algebra embedding $E \hookrightarrow D$ and let $U_E=E^\times \cap U_D^1$. For a character $\psi: E^\times \rightarrow k^\times$,
let $\psi^\dagger$ denote the character of $E^\times U_D^1 \rightarrow k^\times$
such that $\restr{\psi^\dagger}{E^\times} = \psi$ and $\restr{\psi^\dagger}{U_D^1} = 1$
(note that $\psi^\dagger$ is well defined because $U_E$
is a pro-$p$ group, so that $\psi|_{U_E}=1$).
Note that since the abelianization of $D^\times$ is isomorphic to $F^\times$ \cite{matsushima-nakayama},
every character of $D^\times$
factors through the reduced norm map
$\nm_{D^\times/F^\times}: D^\times \rightarrow F^\times$.

\begin{lemma}[{\cite[Proposition 9(a), p. 260]{vigneras-1989-quaternion}, cf. \cite[Proposition 2.5]{tokimoto}}]
    \label{lem:qda-irreps}
    If $\pi$ is an irreducible smooth representation
    of $D^\times$, then either
    \begin{itemize}
    \item $\pi$ is a $1$-dimensional representation
        isomorphic to $\chi \circ \nm$ for
        a character $\chi:F^\times \rightarrow k^\times$.
    \item $\pi$ is a $2$-dimensional representation
        isomorphic to
        $\Ind_{E^\times U_D^1}^{D^\times}(\psi^\dagger)$
        for a character $\psi: E^\times \rightarrow k^\times$
        that does not factor through the norm map
        $\nm_{E/F}: E^\times \rightarrow F^\times$. Moreover, two such representations $\Ind_{E^\times U_D^1}^{D^\times}(\psi^\dagger)$ and $\Ind_{E^\times U_D^1}^{D^\times}(\psi'^\dagger)$ are isomorphic if and only if there is $\sigma \in \Gal(E/F)$ such that $\psi'=\psi^\sigma$.
    \end{itemize}
\end{lemma}

\begin{proof}[Proof of Corollary \ref{cor:qda}]
    We will prove the following stronger statement: for every character $\epsilon$ of $D^\times$, the space $\Hom_{D^\times} \paren{\pi_1 \otimes \pi_2 \otimes \pi_3, \epsilon}$ has dimension at most one. Since, as discussed above, each $\pi_i$ has a twist that is inflated from $G=D^\times/\varpi_F^\Z U_D^1=C \rtimes \abrac{\sigma}$, and since this statement is invariant under twisting, we may assume that each $\pi_i$ is inflated from $G$. Moreover, if at least two of $\pi_1$, $\pi_2$ and $\pi_3$ are characters, then $\pi_{1} \otimes \pi_{2} \otimes \pi_{3}$ is irreducible and there is nothing to show, so we assume, without loss of generality, that both $\pi_1$ and $\pi_2$ have dimension greater than one.  We are now working in the semisimple category of $k[G]$-modules (since $p \nmid \#G)$, and we have to show that the multiplicity of the representation $\pi_{3}^\vee \otimes \epsilon$ in $\pi_{1} \otimes \pi_{2}$ is at most one, where $\pi_{3}^\vee$ is the contragredient of $\pi_{3}$. In fact, $\pi_{1} \otimes \pi_{2}$ is multiplicity-free, as we now show.

    By Lemma \ref{lem:qda-irreps}, there are characters $\psi_1,\psi_2:E^\times \to k^\times$ with $\psi_i \ne \psi_i^\sigma$ such that $\pi_{i}=\Ind_{E^\times U_D^1}^{D^\times}(\psi_i^\dagger)$. Identifying each $\pi_i$ with the associated representation of $G$, this is simply $\pi_i=\Ind_C^G (\psi_i)$. Then
    \begin{align*}
        \pi_1 \otimes \pi_2
            &\cong \Ind_{C}^{G}(\psi_1) \otimes \pi_2 \\
            &\cong \Ind_{C}^{G}\big(\psi_1 \otimes
                \Res_{C}^{G}(\pi_2)\big) \\
            &\cong \Ind_{C}^{G}\big(\psi_1 \otimes
                (\psi_2 \oplus {\psi_2^\sigma})\big) \\
            &\cong \Ind_{C}^{G}(\psi_1 \psi_2)
                \oplus \Ind_{C}^{G}(\psi_1 {\psi_2^\sigma})
        \end{align*}
    by the push-pull formula and Mackey's formula. Since $\psi_i \ne \psi_i^\sigma$, it follows that $\psi_1\psi_2 \ne \psi_1\psi_2^\sigma$ and $\psi_1^\sigma\psi_2^\sigma \ne \psi_1\psi_2^\sigma$, so (using Lemma \ref{lem:qda-irreps} again) this is a direct sum of two distinct two-dimensional representations that do not share a common subrepresentation. In particular, $\pi_1 \otimes \pi_2$ is multiplicity-free.
\end{proof}

%%%%%%%%%%%%%%%%%%%%%%%%%%%% References %%%%%%%%%%%%%%%%%%%%%%%%%%%%%%

\newcommand{\etalchar}[1]{$^{#1}$}

\end{document}